\documentclass[11pt]{article}
\usepackage{mathtools}
\mathtoolsset{centercolon}
\usepackage{microtype}
\usepackage[usenames,dvipsnames,svgnames,table]{xcolor}
\usepackage{amsfonts,amsmath, amssymb,amsthm,amscd}
\usepackage[height=9in,width=6.5in]{geometry}
\usepackage{verbatim}

\usepackage[margin=10pt,font=small,labelfont=bf]{caption}
\usepackage{latexsym}
\usepackage[pdfauthor={Tang},bookmarksnumbered,hyperfigures,colorlinks=true,citecolor=BrickRed,linkcolor=BrickRed,urlcolor=BrickRed,pdfstartview=FitH]{hyperref}
\usepackage{mathrsfs}
\usepackage{graphicx}
\usepackage{enumitem}

\usepackage[inline,nolabel]{showlabels}

\newtheorem{theorem}{Theorem}[section]

\newtheorem{proposition}[theorem]{Proposition}

\newtheorem{example}[theorem]{Example}

\newtheorem{definition}[theorem]{Definition}
\newtheorem{remark}[theorem]{Remark}
\numberwithin{equation}{section}

\newcommand{\be}{\begin{equation}}
	\newcommand{\ee}{\end{equation}}
\newcommand\ba{\begin{align}}
	\newcommand\ea{\end{align}}

\usepackage{mathrsfs}
\usepackage[capitalize]{cleveref}
\Crefname{fact}{Fact}{Facts}
\Crefname{claim}{Claim}{Claims}
\Crefname{example}{Example}{Examples}

\usepackage{tikz}
\usetikzlibrary{decorations.markings}
\tikzset{->-/.style={decoration={
			markings,
			mark=at position .6 with {\arrow{>}}},postaction={decorate}}}
\usepackage{tkz-euclide}
\usetikzlibrary{arrows.meta}

\newtheorem{observation}[theorem]{Observation}

\DeclareMathOperator{\br}{\mathrm{br}}
\DeclareMathOperator{\gr}{\underline{\mathrm{gr}}}
\DeclareMathOperator{\ugr}{\overline{\mathrm{gr}}}

\DeclareMathOperator{\brr}{\mathrm{brr}}
\DeclareMathOperator{\grr}{\underline{\mathrm{grr}}}
\DeclareMathOperator{\ugrr}{\overline{\mathrm{grr}}}

\DeclareMathOperator{\Ibn}{\mathrm{Ibr}}
\DeclareMathOperator{\Igr}{\underline{\mathrm{Igr}}}
\DeclareMathOperator{\uIgr}{\overline{\mathrm{Igr}}}

\begin {document}
\title{A subperiodic tree whose intermediate branching number \\ is strictly less than the lower intermediate growth rate}
\author{
	Pengfei Tang\thanks{Department of Mathematical Sciences, Tel Aviv University. This research is supported by ERC consolidator grant 101001124 (UniversalMap), and by ISF grant 1294/19.
		Current affiliation: Center for Applied Mathematics, Tianjin University, China. 
		Email: \textsf{pengfei\_tang@tju.edu.cn}.}
}
\date{}
\maketitle

% MSC: 60J10  	Markov chains (discrete-time Markov processes on discrete state spaces)
% 05C81  	Random walks on graphs
% 05C05  	Trees

% keywords: trees; cover time; BEST theorem

\begin{abstract}
We construct an example of a subperiodic tree whose intermediate branching number is strictly less than the lower intermediate growth rate. This answers a question of Amir and Yang (2022) in the negative.
\end{abstract}

\section{Introduction and main result}

There are several ways to measure the branching structure of an infinite locally finite tree. An important and successful one is the branching number introduced by Lyons \cite{Lyons1990rwpt}. For instance the branching number is the critical parameter for Bernoulli percolation and homesick random walk on trees. 
However the branching number is not so effective for trees with sub-exponential growth. Later Collevecchio,  Kious and Sidoravicius \cite{Collevecchio_etal2020branching_ruin_number} introduced a branching-ruin number which works well for trees with polynomial  growth. Inspired by these previous work,  recently Amir and Yang \cite{Amir_Yang2022branching} introduced the intermediate branching number and showed that it is crucial for several probability models on trees with intermediate growth rate. 

Our focus here is a special family of infinite locally finite trees---the subperiodic trees. For a subperiodic tree, the branching number actually equals the exponential growth rate---this result is due to Furstenberg \cite{Furstenberg1967}; see Theorem 3.8 in \cite{LP2016} for a proof. Amir and Yang \cite{Amir_Yang2022branching} then asked whether the corresponding equality holds for the intermediate branching number and the lower intermediate growth rate on subperiodic trees. In the present note we construct an example of a subperiodic tree whose intermediate branching number is strictly less than its lower intermediate growth rate, answering their question in the negative.

\subsection{Various branching numbers and growth rates of infinite trees}
Suppose $T=(V,E)$ is an infinite locally finite tree with a distinguished vertex $o$, which will be called the \textbf{root} of $T$.  We imagine the tree $T$ as growing upward from the root $o$. For $x,y\in V$, we write $x\leq y$ if $x$ is on the shortest path from $o$ to $y$; and $T^x$ for the \textbf{subtree} of $T$ containing all the vertices $y$ with $y\geq x$. For a vertex $x\in V$ we denote by $|x|$ the graph distance from $o$ to $x$. For an edge $e\in E$, we write $e=(e^-,e^+)$ where $|e^+|=|e^-|+1$ and define $|e|=|e^+|$. Write $T_n:=\{e\in E\colon |e|=n \}$. Write $B(n)=\{x\colon x\in V, |x|\leq n\}$ for the ball of radius $n$ centered at $o$.

A \textbf{cutset} $\pi$ separating $o$ and infinity is a set of edges such that every infinite path starting from $o$ must include an edge in $\pi$. For instance $T_n$ is a cutset separating $o$ and infinity for every $n\geq1$.  We write $\Pi(T)$ for the collection of cutsets separating $o$ and infinity. 
The \textbf{branching number} of $T$ is defined as 
\be\label{eq: def of branching number}
\br(T):=\sup\left\{ \lambda>0\colon \inf_{\pi\in\Pi(T)}\sum_{e\in\pi} \lambda^{-|e|}>0 \right\}.
\ee
We recommend the readers Chapter 3 of \cite{LP2016} for backgrounds on branching numbers. The \textbf{lower exponential growth rate} of $T$ is defined as 
\be\label{eq: def of lower exponential growth rate}
\gr(T):=\liminf_{n\to\infty}|T_n|^{1/n}.
\ee
The \textbf{upper exponential growth rate} of $T$ is defined as $\ugr(T):=\limsup_{n\to\infty}|T_n|^{1/n}$ similarly.
Note that $\gr(T)$ can  be rewritten in a similar form as \eqref{eq: def of branching number}:
\[
\gr(T)=\sup\left\{ \lambda>0\colon \liminf_{n\to\infty}\sum_{e\in T_n}\lambda^{-|e|}>0 \right\}
\]
and in particular
\[
1\leq \br(T)\leq \gr(T). 
\]

The \textbf{branching-ruin number} introduced by Collevecchio,  Kious and Sidoravicius \cite{Collevecchio_etal2020branching_ruin_number} is defined as 
\be\label{eq: def branching-ruin}
\brr(T):=\sup\left\{ \lambda>0\colon \inf_{\pi\in\Pi(T)}\sum_{e\in\pi} |e|^{-\lambda}>0 \right\},
\ee
where we use the convention of $\sup \emptyset=0$. This branching-ruin number is a natural way to measure trees with polynomial growth rate and turned out be the critical parameter of some random processes \cite{Collevecchio_etal2019EJP} (in particular  the once-reinforced random walk \cite{Collevecchio_etal2020branching_ruin_number}). One can define  corresponding \textbf{lower (upper) polynomial growth rates} by 
\be\label{eq: def of growth-ruin}
\grr(T):=\liminf_{n\to\infty}\frac{\log |T_n|}{\log n} \quad \textnormal{ and }\quad \ugrr(T):=\limsup_{n\to\infty}\frac{\log |T_n|}{\log n}.
\ee
Note that
\[
\grr(T)=\sup\left\{ \lambda>0\colon \liminf_{n\to\infty}\sum_{e\in T_n}|e|^{-\lambda}>0 \right\}.
\]
and in particular $\brr(T)\leq \grr(T)$.

Recently Amir and Yang \cite{Amir_Yang2022branching} introduced the \textbf{intermediate branching number} 
\be\label{eq: def intermediate branching}
\Ibn(T):=\sup\left\{ \lambda>0\colon \inf_{\pi\in\Pi(T)}\sum_{e\in\pi} \exp\big(-|e|^\lambda\big)>0 \right\}
\ee
and the \textbf{lower (upper)  intermediate growth rates}
\be\label{eq: def intermediate growth}
\Igr(T):=\liminf_{n\to\infty}\frac{\log\log |T_n|}{\log n} \quad \textnormal{ and }\quad \uIgr(T):=\limsup_{n\to\infty}\frac{\log\log |T_n|}{\log n}.
\ee
Similarly,  
\[
\Ibn(T)\leq \Igr(T)=\sup\left\{ \lambda>0\colon \liminf_{n\to\infty}\sum_{e\in T_n}\exp\big(-|e|^\lambda\big)>0 \right\}.
\]
Amir and Yang \cite{Amir_Yang2022branching} proved that the intermediate branching number is the critical parameter for certain random walk, percolation and firefighting problems on trees with intermediate growth, where a tree $T$ was said to be of \textbf{intermediate (stretched exponential) growth} if $0<\Igr(T)\leq \uIgr(T)<1.$

We remark that these numbers $\br(T),\gr(T),\ugr(T),\brr(T),\grr(T),\ugrr(T),\Ibn(T),\Igr(T)$ and $\uIgr(T)$ do not depend on the choice of the root of $T$. 

\subsection{Subperiodic trees}

We first recall the definition of subperiodic trees from p~82 of \cite{LP2016}; see Example 3.6 and 3.7 there for some examples of subperiodic trees. 
\begin{definition}\label{def: subperiodic trees}
	Let $N\in \{0,1,2,3,\ldots\}$. An infinite tree $T$ is called \textbf{$N$-periodic} (resp., \textbf{$N$-subperiodic}) if $\forall \,x\in T$ there exists an adjacency-preserving bijection (resp.\, injection) $f: T^x\to T^{f(x)}$ with $|f(x)|\leq N$. A tree is \textbf{periodic} (resp.\, \textbf{subperiodic}) if there is some $N$ for which it is $N$-periodic (resp., $N$-subperiodic). 
\end{definition}

As mentioned earlier  $\br(T)=\gr(T)=\ugr(T)$ for any subperiodic tree $T$ (\cite[Theorem 3.8]{LP2016}). Amir and Yang noticed that there exist subperiodic trees such that $\Igr(T)<\uIgr(T)$ (see Section 4.1 of \cite{Amir_Yang2022branching}) and asked\footnote{See (2.12) on page 4 of version 3 of the paper \cite{Amir_Yang2022branching} on arxiv}  whether $\Ibn(T)=\Igr(T)$ for subperiodic trees with intermediate growth rate. Our main result gives a negative answer to their question.
\begin{theorem}\label{thm: an example with Ibn less than Igr}
	There exists a subperiodic tree $T$ with intermediate growth rate and 
	\[
	\Ibn(T)<\Igr(T).
	\]
\end{theorem}

\section{Proof of the main result}

We will prove  \cref{thm: an example with Ibn less than Igr} via a concrete example (see Example \ref{exam: Ibn less than Igr}).

\subsection{Coding by trees}

%
%Given an integer $b\geq2$, one can code a closed nonempty set $E\subset [0,1]$ as a subtree $T=T_{[b]}(E)$ of the \textbf{$b$-ary tree} $\mathbb{T}_b$ in the following way. The vertices of $T$ correspond to the set of $b$-adic intervals with nonempty intersections with $E$, where an interval of the form $[k/b^n,(k+1)/b^n]$ for integers $k$ and $n$ is called a \textbf{$b$-adic interval of order $n$}. We let the root of $T$ be the vertex corresponding to $[0,1]$. Two such intervals are connected by an edge if and only if one interval contains the other and the orders of them differ by $1$. See Section 1.10 and 15.2 of \cite{LP2016} for background on this coding. 
%
%Note that if a point $x\in E$ has the form of $k/b^n$ (i.e., it is the endpoint of some $b$-adic interval), then it might correspond to two rays in $T_{[b]}(E)$. This is just the fact that  $x$ can be written in base $b$ with two equivalent expressions $x=\frac{k}{b^n}=\frac{k-1}{b^n}+\sum_{m=n+1}^{\infty}\frac{b-1}{b^m}$. 
%To build a bijection between the rays of $T$ and $E$ we will view $T$ as a labelled subtree of $\mathbb{T}_b$ and $E$ as a subset of $\{0,1,\ldots,b-1\}^{\mathbb{N}}$ instead, where $\mathbb{N}=\{1,2,3,\ldots\}$.

%We will only consider the $b=3$ case for simplicity and view $\mathbb{T}_3$ as a labelled tree with the root labelled as $\emptyset$, the three children of the root labelled $0,1,2$ respectively from left to right, and so on.

Our example will be a subtree of the \textbf{$3$-ary tree} $\mathbb{T}_3$ and we view $\mathbb{T}_3$ as a \textbf{labelled} tree with the root labelled as $\emptyset$, the three children of the root labelled $0,1,2$ respectively from left to right, and so on. Write $\mathscr{D}(\mathbb{T}_3)$ for the set of infinite labelled subtrees of $\mathbb{T}_3$ which contain the root and have no leaf and write $\mathscr{R}(\mathbb{T}_3)$ for the set of labelled rays starting from the root. In particular $\mathscr{R}(\mathbb{T}_3)\subset \mathscr{D}(\mathbb{T}_3)$.

For each element $a=(a_1,a_2,a_3,\ldots)\in\{0,1,2\}^{\mathbb{N}}$, we associate it with a ray $\Phi(a)\in\mathscr{R}(\mathbb{T}_3)$ with the $(n+1)$-th vertex on the ray labelled as $a_1a_2\cdots a_n$. (The first vertex is just the root labelled as $\emptyset$.)  Obviously $\Phi$ is a bijection between $\{0,1,2\}^{\mathbb{N}}$ and $\mathscr{R}(\mathbb{T}_3)$. We now extend $\Phi$ as a mapping from all nonempty subsets of $\{0,1,2\}^{\mathbb{N}}$ to $\mathscr{D}(\mathbb{T}_3)$.
\begin{definition}[Coding by trees]\label{def: map from E to T}
	For a  nonempty subset $E$ of $\{0,1,2\}^{\mathbb{N}}$, the tree $\Phi(E)\in\mathscr{D}(\mathbb{T}_3)$ is defined as the union of the rays (each ray is viewed as a labelled subtree of $\mathbb{T}_3$) $\Phi(E)=\bigcup_{x\in E}\Phi(x)$, where the union means the vertex set of $\Phi(E)$ is the union of the vertex set of $\Phi(x)$ and the same for the edge set. 
\end{definition}
\begin{remark}
	One can define a natural metric $d$ on $\{0,1,2\}^{\mathbb{N}}$: the distance between two elements $x=(x_1,x_2,\ldots)$ and $y=(y_1,y_2,\ldots)\in \{0,1,2\}^{\mathbb{N}}$ is given by
	\[
	d(x,y):=\frac{1}{e^k}, \quad \textnormal{ where }k=k(x,y)=\inf\{ i\colon x_i\neq y_i\}.
	\]

	The coding by trees in Definition \ref{def: map from E to T} comes from the canonical coding of closed subsets of the interval $[0,1]$ by trees; see Section 1.10 and 15.2 of \cite{LP2016} for background. We simply replace $[0,1]$ by the metric space $\big( \{0,1,2\}^{\mathbb{N}},d~\big)$ for convenience. 
\end{remark}

It is straightforward to verify that 
\begin{itemize}
	\item For each nonempty subset $E\subset \{0,1,2\}^{\mathbb{N}}$ and its closure $\overline{E}$ in the metric space $\big( \{0,1,2\}^{\mathbb{N}},d~\big)$, one has that $\Phi(E)=\Phi\big(\overline{E}\big)$.
	
	\item Moreover the map $\Phi$ is a bijection if its domain is restricted to the collection of all nonempty  \textbf{closed} subsets of $\{0,1,2\}^{\mathbb{N}}$.
\end{itemize}

We also define the \textbf{shift map} $\mathcal{S}:\{0,1,2\}^{\mathbb{N}} \to \{0,1,2\}^{\mathbb{N}}$ by 
\[
\mathcal{S}\big((a_1,a_2,a_3,\ldots)\big)=(a_2,a_3,a_4,\ldots). 
\]

The following observation is a rephrasing of Example 3.7 in \cite{LP2016} in the case $b=3$ and it is crucial for our construction later. 
\begin{observation}\label{obser: invariance imply subperiodic}
	If a nonempty closed subset $E\subset \{0,1,2\}^{\mathbb{N}}$ is invariant under the shift map in the sense that $\mathcal{S}(E)\subset E$, then the tree $\Phi(E)$ is $0$-subperiodic. 
\end{observation}

\subsection{The construction of our example}
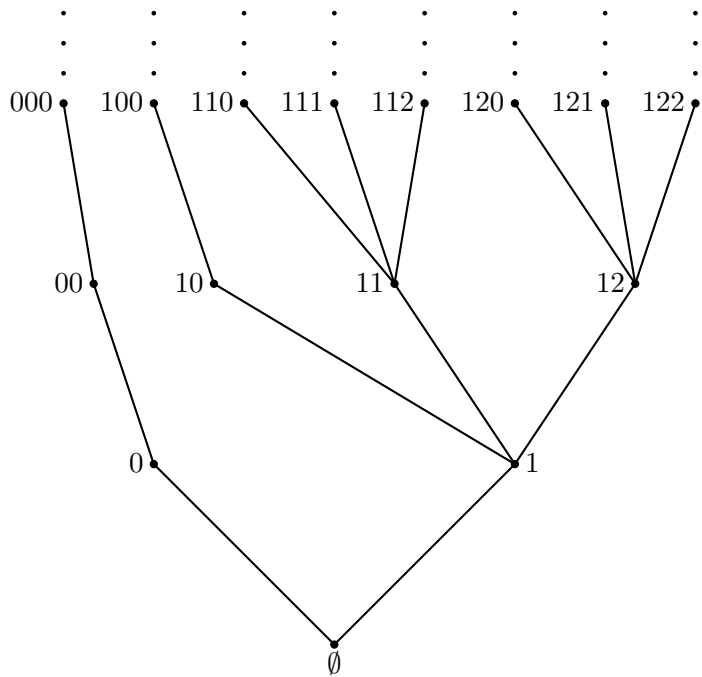
\begin{figure}[h!]
	\centering
	\begin{tikzpicture}[scale=0.8, text height=1.5ex,text depth=.25ex] 
		%\draw [help lines] (0,0) grid (14,12);

		%%%%%%%%%% the root %%%%%%%%%%%
		\draw[fill=black] (6,0) circle [radius=0.06];
		\node[below] at (6,0) {$\emptyset$};
		
		%%%%%%%%%%%%%  T_1 %%%%%%%%%%%%%
		\draw[fill=black] (3,3) circle [radius=0.06];
		\node[left] at (3,3) {$0$};

		\draw[fill=black] (9,3) circle [radius=0.06];
		\node[right] at (9,3) {$1$};
		
		\draw[color=black,thick] (6,0)--(3,3);
		\draw[color=black,thick] (6,0)--(9,3);

		%%%%%%%%%%%%%%% T_2 %%%%%%%%%%%%%%
		
		\draw[fill=black] (2,6) circle [radius=0.06];
		\node[left] at (2,6) {$00$};

		\draw[fill=black] (4,6) circle [radius=0.06];
		\node[left] at (4,6) {$10$};
		
		\draw[fill=black] (7,6) circle [radius=0.06];
		\node[left] at (7,6) {$11$};
		
		\draw[fill=black] (11,6) circle [radius=0.06];
		\node[left] at (11,6) {$12$};
		
		\draw[color=black,thick] (2,6)--(3,3);
		\draw[color=black,thick] (9,3)--(4,6);
		\draw[color=black,thick] (9,3)--(7,6);
		\draw[color=black,thick] (9,3)--(11,6);

		%%%%%%%%%%%%%%%%% T_3 %%%%%%%%%%%%%%%%%%
		\draw[fill=black] (1.5,9) circle [radius=0.06];
		\node[left] at (1.5,9) {$000$};
		\draw[color=black,thick] (2,6)--(1.5,9);
		
		\draw[fill=black] (3,9) circle [radius=0.06];
		\node[left] at (3,9) {$100$};
		\draw[color=black,thick] (4,6)--(3,9);
		
		\draw[fill=black] (4.5,9) circle [radius=0.06];
		\node[left] at (4.5,9) {$110$};
		\draw[color=black,thick] (7,6)--(4.5,9);
		
		\draw[fill=black] (6,9) circle [radius=0.06];
		\node[left] at (6,9) {$111$};
		\draw[color=black,thick] (7,6)--(6,9);
		
		\draw[fill=black] (7.5,9) circle [radius=0.06];
		\node[left] at (7.5,9) {$112$};
		\draw[color=black,thick] (7,6)--(7.5,9);
		
		\draw[fill=black] (9,9) circle [radius=0.06];
		\node[left] at (9,9) {$120$};
		\draw[color=black,thick] (11,6)--(9,9);
		
		\draw[fill=black] (10.5,9) circle [radius=0.06];
		\node[left] at (10.5,9) {$121$};
		\draw[color=black,thick] (11,6)--(10.5,9);
		
		\draw[fill=black] (12,9) circle [radius=0.06];
		\node[left] at (12,9) {$122$};
		\draw[color=black,thick] (11,6)--(12,9);

		\foreach \x in {1.5, 3, 4.5, 6, 7.5, 9, 10.5, 12}
		\foreach \y in {9.5, 10, 10.5}
		{
			\draw[fill=black] (\x,\y) circle [radius=0.03];
		}
		
	\end{tikzpicture}
	\caption{The $1$-$3$ tree $T_{1,3}$ and its labeling.}
	\label{fig: one-three tree}
\end{figure}

We first review the $1$-$3$ tree $T_{1,3}$  \cite[Example 1.2]{LP2016}: the root has two children; and $|T_n|=2^n$; and for each $n\geq 1$, the left half vertices at distance $n$ from the root will each have only $1$ child, the right half will each have $3$ children. We view $T_{1,3}$ as a labelled subtree of $\mathbb{T}_3$ according to the following \textbf{labeling rule}: the root is labelled as $\emptyset$ and if  a vertex with label $a_1a_2\cdots a_n$ has $k$ children, then its $k$ children are labelled as $a_1a_2\cdots a_n0,\ldots,a_1a_2\cdots a_n(k-1)$ respectively from left to right.  See Figure \ref{fig: one-three tree} for $T_{1,3}$ and its labeling. 

\begin{figure}[h!]
	\centering
	\begin{tikzpicture}[scale=0.8, text height=1.5ex,text depth=.25ex] 
		%\draw [help lines] (0,0) grid (14,12);

		%%%%%%%%%% the root %%%%%%%%%%%
		\draw[fill=black] (6,0) circle [radius=0.06];
		\node[below] at (6,0) {$\emptyset$};
		
		%%%%%%%%%%%%% old T_1 %%%%%%%%%%%%%
		\draw[fill=black] (3,3) circle [radius=0.06];
		\node[left] at (3,3) {$0$};

		\draw[fill=black] (9,3) circle [radius=0.06];
		\node[right] at (9,3) {$1$};
		
		\draw[color=black,thick] (6,0)--(3,3);
		\draw[color=black,thick] (6,0)--(9,3);

		%%%%%%%%%%%% new T_2 %%%%%%%%%%%
		\draw[fill=black] (2.5,4.5) circle [radius=0.06];
		\node[left] at (2.5,4.5) {$00$};
		
		\draw[fill=black] (6.5,4.5) circle [radius=0.06];
		\node[left] at (6.5,4.5) {$10$};
		
		\draw[fill=black] (8,4.5) circle [radius=0.06];
		\node[left] at (8,4.5) {$11$};
		
		\draw[fill=black] (10,4.5) circle [radius=0.06];
		\node[left] at (10,4.5) {$12$};

		%%%%%%%%%%%%%%% old T_2= new T_3  %%%%%%%%%%%%%%
		
		\draw[fill=black] (2,6) circle [radius=0.06];
		\node[left] at (2,6) {$000$};

		\draw[fill=black] (4,6) circle [radius=0.06];
		\node[left] at (4,6) {$100$};
		
		\draw[fill=black] (7,6) circle [radius=0.06];
		\node[left] at (7,6) {$110$};
		
		\draw[fill=black] (11,6) circle [radius=0.06];
		\node[right] at (11,6) {$120$};
		
		\draw[color=black,thick] (2,6)--(3,3);
		\draw[color=black,thick] (9,3)--(4,6);
		\draw[color=black,thick] (9,3)--(7,6);
		\draw[color=black,thick] (9,3)--(11,6);

		%%%%%%%%%%%% new T_4 and T_5 %%%%%%%%%

		\draw[fill=black] (2-1/6,6+1) circle [radius=0.06];
		\node[left] at (2-1/6,6+1) {$0000$};
		\draw[fill=black] (2-1/3,6+2) circle [radius=0.06];
		\node[left] at (2-1/3,6+2)  {$00000$};

		\draw[fill=black] (4-1/3,6+1) circle [radius=0.06];
		\node[left] at (4-1/3,6+1) {$1000$};
		\draw[fill=black] (4-2/3,6+2) circle [radius=0.06];

		\draw[fill=black] (7-2.5/3,6+1) circle [radius=0.06];
		\node[left] at (7-2.5/3,6+1) {$1100$};
		\draw[fill=black] (7-5/3,6+2) circle [radius=0.06];

		\draw[fill=black] (7-1/3,6+1) circle [radius=0.06];
		%\node[above] at (7-1/3,6+1) {$1101$};
		\draw[fill=black] (7-2/3,6+2) circle [radius=0.06];

		\draw[fill=black] (7+0.5/3,6+1) circle [radius=0.06];
		\node[right] at (7+0.5/3,6+1) {$1102$};
		\draw[fill=black] (7+1/3,6+2) circle [radius=0.06];

		\draw[fill=black] (11-2/3,6+1) circle [radius=0.06];
		\node[left] at (11-2/3,6+1) {$1200$};
		\draw[fill=black] (11-4/3,6+2) circle [radius=0.06];

		\draw[fill=black] (11-0.5/3,6+1) circle [radius=0.06];
		%\node[above] at (11-0.5/3,6+1) {$1201$};
		\draw[fill=black] (11-1/3,6+2) circle [radius=0.06];
		
		\draw[fill=black] (11+1/3,6+1) circle [radius=0.06];
		\node[right] at (11+1/3,6+1) {$1202$};
		\draw[fill=black] (11+2/3,6+2) circle [radius=0.06];
		\node[right] at (11+2/3,6+2) {$12020$};

		%%%%%%%%%%%%%%%%% old T_3 =new T_6  %%%%%%%%%%%%%%%%%%
		\draw[fill=black] (1.5,9) circle [radius=0.06];
		\node[left] at (1.5,9) {$000000$};
		\draw[color=black,thick] (2,6)--(1.5,9);
		
		\draw[fill=black] (3,9) circle [radius=0.06];
		%\node[left] at (3,9) {$100000$};
		\draw[color=black,thick] (4,6)--(3,9);
		
		\draw[fill=black] (4.5,9) circle [radius=0.06];
		%\node[left] at (4.5,9) {$110$};
		\draw[color=black,thick] (7,6)--(4.5,9);
		
		\draw[fill=black] (6,9) circle [radius=0.06];
		%\node[left] at (6,9) {$111$};
		\draw[color=black,thick] (7,6)--(6,9);
		
		\draw[fill=black] (7.5,9) circle [radius=0.06];
		%\node[left] at (7.5,9) {$112$};
		\draw[color=black,thick] (7,6)--(7.5,9);
		
		\draw[fill=black] (9,9) circle [radius=0.06];
		%\node[left] at (9,9) {$120$};
		\draw[color=black,thick] (11,6)--(9,9);
		
		\draw[fill=black] (10.5,9) circle [radius=0.06];
		%\node[left] at (10.5,9) {$121$};
		\draw[color=black,thick] (11,6)--(10.5,9);
		
		\draw[fill=black] (12,9) circle [radius=0.06];
		\node[right] at (12,9) {$120200$};
		\draw[color=black,thick] (11,6)--(12,9);

		\foreach \x in {1.5, 3, 4.5, 6, 7.5, 9, 10.5, 12}
		\foreach \y in {9.5, 10, 10.5}
		{
			\draw[fill=black] (\x,\y) circle [radius=0.03];
		}
		
	\end{tikzpicture}
	\caption{The tree $T_{0}$ and its labeling.}
	\label{fig: T zero}
\end{figure}
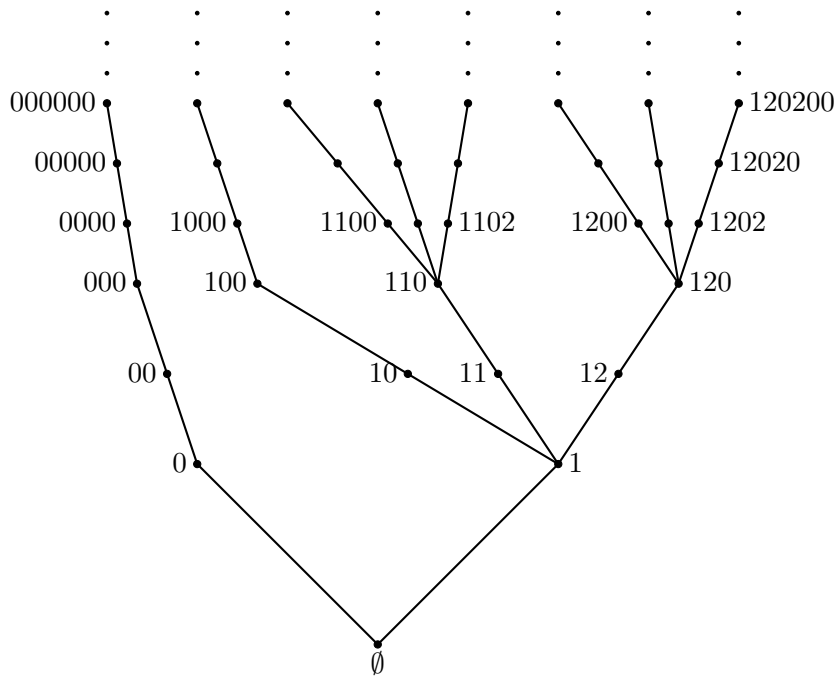

\begin{example}\label{exam: Ibn less than Igr}
	Let $T_{0}$ be the tree obtained by replacing each edge $e$ of the $1$-$3$ tree $T_{1,3}$  by a path of length $|e|$ and view it as a  subtree of $\mathbb{T}_3$ labelled according to the labeling rule we used for $T_{1,3}$ (see Figure \ref{fig: T zero}). 	As already noted by Amir and Yang \cite{Amir_Yang2022branching}, the tree $T_0$ satisfies
	\be\label{eq: Ibn and Igr of Tzero}
	\Ibn(T_0)=0\quad \textnormal{ and }\quad \Igr(T_0)=\frac{1}{2}.
	\ee
	However $T_0$ is not subperiodic. 
	
	Let $E_0$ be the closed subset of $\{0,1,2\}^{\mathbb{N}}$ such that $\Phi(E)=T_0$. 
	Define $E_j=\mathcal{S}(E_{j-1})$ for $j\geq1$ and let $\widetilde{E}=\bigcup_{j=0}^{\infty}E_j$. 
	Our example is just the tree $\widetilde{T}:=\Phi\big(\widetilde{E}\big)$. 
\end{example}

Recall that $\mathscr{D}(\mathbb{T}_3)$ denotes the set of infinite labelled subtrees of $\mathbb{T}_3$ which contain the root and have no leaf.
For a vertex $v\in V(T_0)$ labelled as $a_1a_2\cdots a_n$, we will view the subtree $T_0^v$ as  a labelled subtree of $\mathbb{T}_3$ rooted at $\emptyset$, i.e., view it as the tree
\[
\Phi\big(\mathcal{S}^{\circ n}\{ x=(x_1,x_2,x_3,\ldots)\colon x\in E_0,x_i=a_i \textnormal{ for } i=1,\ldots,n  \}   \big) \in \mathscr{D}(\mathbb{T}_3).
\]
Since   $E_n=\bigcup_{v\in V(T_0),|v|=n}\mathcal{S}^{\circ n}\{ x=(x_1,x_2,x_3,\ldots)\colon x\in E_0,v \textnormal{ is labelled as }x_1x_2\cdots x_n  \}$, 
we have the following equivalent description of $\widetilde{T}$: 
\begin{observation}\label{obser: T tilde as union}
	As labelled subtrees of $\mathbb{T}_3$, the tree $\widetilde{T}$ is just the union of $T_0^v$ over all $v\in V(T_0)$. 
\end{observation}

\subsection{The intermediate branching number and the intermediate growth rate of our example}
By construction the set $\widetilde{E}$ is invariant under the shift map. Thus by \cref{obser: invariance imply subperiodic}  the tree $\widetilde{T}=\Phi\big(\widetilde{E}\big)$ is subperiodic. We will show that $0=\Ibn\big(\widetilde{T}\big)< \Igr\big(\widetilde{T}\big)=\frac{1}{2}$
which then proves \cref{thm: an example with Ibn less than Igr}.

\begin{proposition}
	For the tree $\widetilde{T}=\Phi\big(\widetilde{E}\big)$ constructed in Example \ref{exam: Ibn less than Igr}, one has that 
	\[
	\Igr\big(\widetilde{T}\big)=\uIgr\big(\widetilde{T}\big)=\frac{1}{2}\quad \textnormal{ and }\quad \Ibn\big(\widetilde{T}\big)=0.
	\]
\end{proposition}
\begin{proof}
	First of all since $T_0=\Phi(E_0)$ is a subtree of $\widetilde{T}=\Phi\big(\widetilde{E}\big)$, one has that
	\[
	\uIgr\big(\widetilde{T}\big)\ge \Igr\big(\widetilde{T}\big)\geq \Igr(T_0)\stackrel{\eqref{eq: Ibn and Igr of Tzero}}{=}\frac{1}{2}.
	\]
	On the other hand, note that $	\big|\widetilde{T}_n\big|$ equals the cardinality of the set  $\{(x_1,\ldots,x_n)\colon x=(x_1,x_2,\ldots)\in\widetilde{E} \}$---the first $n$-bits of $\widetilde{E}$.
	Also observe that
	a ray $\gamma$ in $T_{1,3}$ coding the sequence $(a_1,a_2,a_3,\ldots)$ becomes a ray $\gamma'$ in $T_0$ coding the sequence 
	\[
	(a_1,a_2,0,a_3,0,0,a_4,0,0,0,a_5,0,0,0,0,a_6,0,\ldots). 
	\]
	Hence by our construction an element $a\in\widetilde{E}$ always has the form 
	\be\label{eq: general form of an element in E zero}
	a=\big(\underbrace{0,\ldots,0}_{m},a_j,\underbrace{0,\ldots,0}_{= j-1},a_{j+1},\underbrace{0,\ldots,0}_{= j},a_{j+2},0,\cdots
	\big),
	\ee
	where $(a_1,a_2,a_3,\ldots)\in\Phi^{-1}(T_{1,3})$ and $m\leq \max(j-2,0)$. Note that there exists a constant $c>0$ such that there are at most $c\sqrt{n}+1$  nontrivial entries $a_j,a_{j+1},\ldots,a_{j+c\sqrt{n}}$ in the first $n$-bits of $a$. If $j\geq n+1$, then there is at most one nonzero entry in the first $n$-bits and this would contribute at most $2n+1$ to the set $\{(x_1,\ldots,x_n)\colon x=(x_1,x_2,\ldots)\in\widetilde{E} \}$. If $j\leq n$, then there are at most $\max(n-2,0)\leq n$ choices for  $m$---the number of zeroes   before $a_j$; once $m$ and $j$ are fixed, the positions of $a_j,a_{j+1},\ldots,a_{j+c\sqrt{n}}$ are fixed and each element of  $\{a_j,a_{j+1},\ldots,a_{j+c\sqrt{n}}\}$ has at most $3$ choices, hence this contributes at most $n^2*3^{c\sqrt{n}+1}$ to the set $\{(x_1,\ldots,x_n)\colon x=(x_1,x_2,\ldots)\in\widetilde{E} \}$. In sum  we have  $\big|\widetilde{T}_n\big|\leq 3^{C\sqrt{n}}$  for some constant $C>0$. Therefore one has the other direction
	\[
	\Igr\big(\widetilde{T}\big)\leq \uIgr\big(\widetilde{T}\big)\leq  \frac{1}{2}.
	\]
	
	Next we proceed to show that $\Ibn\big(\widetilde{T}\big)=0$. Fixing an arbitrary $\lambda>0$, we will show that for any $\varepsilon>0$ there exists a cutset $\pi$ of $\widetilde{T}$ such that 
	\be\label{eq: capacity arbitrary small}
	\sum_{e\in\pi} \exp\big(-|e|^\lambda\big)\leq 2\varepsilon. 
	\ee
	
	Since $\Ibn(T_0)\stackrel{\eqref{eq: Ibn and Igr of Tzero}}{=}0$, one has $\Ibn(T_0^v)=0$ for any $v\in V(T_0)$. In particular one can choose  cutsets  $\pi_v$ for $T_0^v$ (viewed as a subtree of $\mathbb{T}_3$ rooted at $\emptyset$) such that 
	\[
	\sum_{v \in V(T_0)} \sum_{e\in\pi_v}\exp\big(-|e|^\lambda\big)\leq \varepsilon. 
	\]
	Since $\widetilde{T}$ is the union of $T_0^v$ over $v\in V(T_0)$ (Observation \ref{obser: T tilde as union}), one might hope the set $\bigcup_{v\in V(T_0)}\pi_v$ is a cutset of $\widetilde{T}$. But it might  not  be the case since there might exist a ray $\gamma$ in $\widetilde{T}$ such that its  edges come from $T_0^{v_i}$ for infinitely many different $v_i$'s and $\gamma$ is not blocked from infinity by $\bigcup_{v\in V(T_0)}\pi_v$. 
	To rescue this, we add some additional edges in the following way. 
	Choose $N=N(\lambda,\varepsilon)$ large enough so that $9N\exp(-N^\lambda)\leq\varepsilon$. Let $\beta$ be the collection of all edges in $\widetilde{T}_{N+1}$ with the form 
	\[
	(v,vj)\colon v=(v_1v_2\cdots v_N) \textnormal{ with at most one nonzero entry and }j=0,1,2. 
	\]
	In particular 
	\[
	\sum_{e\in\beta}\exp\big(-|e|^\lambda\big)\leq 9N\exp(-N^\lambda)\leq\varepsilon. 
	\]
	Now we set 
	\[
	\pi=\Big(\bigcup_{v\in V(T_0)}\pi_v\Big)\cup \beta.
	\]
	and claim that $\pi$ is  a cutset of $\widetilde{T}$. In fact since $\widetilde{T}$ is just the union of $T_0^v$ over all $v\in V(T_0)$, we can choose $M\geq 100N^2$ large enough so that all the edges $e$ of $\widetilde{T}$ with $|e|\leq N$ appear in some $T_0^v$ with $|v|\leq M$. Now if a ray $\gamma$ of $\widetilde{T}$ does not use any edge outside $\bigcup_{v\in V(T_0),|v|\leq M}T_0^v$, then there must exist some $v\in V(T_0)$ with $|v|\leq M$ such that $\gamma$ is just a ray in $T_0^v$. Hence in this case $\gamma$  has a nonempty intersection with $\pi_v$. Otherwise $\gamma$ must use some edge $e'$ of $\widetilde{T}$ which is not in the union  $\bigcup_{v\in V(T_0),|v|\leq M}T_0^v$. By our choice of $M$, one must have $|e'|>N$ and $e'$ is coming from some $T_0^v$ with $|v|>M\geq 100N^2$. For such a vertex $v$,  in the first $N$ levels of $T_0^v$ there is at most one vertex with three children because of the long pieces of zeroes (see \eqref{eq: general form of an element in E zero} and Figure \ref{fig: T zero}). Therefore the edge $e'$ must be a descendant of some edge from the set $\beta$ and so $\gamma$ has a  nonempty intersection with $\beta$. Hence $\pi$ is a cutset of $\widetilde{T}$.

	By our choice of $\pi_v$ and $\beta$  the cutset $\pi$ satisfies \eqref{eq: capacity arbitrary small}. 
	By \eqref{eq: capacity arbitrary small} one obtains that $	\Ibn\big(\widetilde{T}\big)\leq \lambda$. Since this is true for any $\lambda>0$ one has that $\Ibn\big(\widetilde{T}\big)=0$. 
\end{proof}

\section{Concluding remarks}
In the construction of $T_0$ we replace an edge $e$ by a path of length $f(|e|)$ where the function $f:\mathbb{N}\to\mathbb{N}$ is given by $f(x)=x$. If we use some other increasing functions, say $f(x)=\lceil x^s\rceil$ with $s\in(0,\infty)$, then we can obtain a family of subperiodic trees using the procedure in Example \ref{exam: Ibn less than Igr} so that for each $\alpha\in(0,1)$ there are some trees $T$ in the family with the property that $0=\Ibn(T)<\Igr(T)=\alpha$.

\begin{figure}[h!]
	\centering
	\begin{tikzpicture}[scale=0.8, text height=1.5ex,text depth=.25ex] 
		%\draw [help lines] (0,0) grid (14,12);
		
		\draw (0,3)--(7,3);
		%\draw [dotted] (0,0)--(0,6);
		%\draw[dotted] (7,0)--(7,6);
		
		\foreach \x in {1, 2, 3, 4, 5, 6}
		{\draw (\x,0)--(\x,6);}
		
		\foreach \x in {1, 2, 3, 4, 5, 6}
		\foreach \y in {1, 2, 3, 4, 5}
		{
			\draw[fill=black] (\x,\y) circle [radius=0.05];
		}
		
	\end{tikzpicture}
	\caption{A lexicographically minimal spanning tree of $\mathbb{Z}^2$.}
	\label{fig: lexico Z two}
\end{figure}
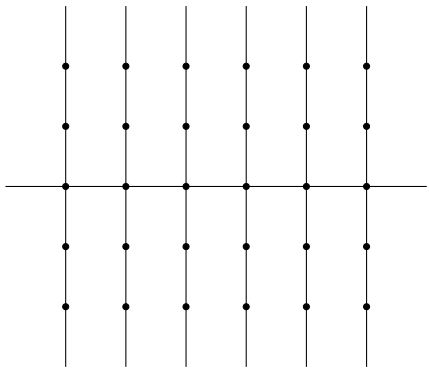

We also note that there exist periodic trees $T$ with polynomial growth that satisfy $\brr(T)<\grr(T)$. For instance consider the following \textbf{lexicographically minimal spanning tree}  of $\mathbb{Z}^2$ illustrated in Figure \ref{fig: lexico Z two}; see Section 3.4 in \cite{LP2016} for definitions of Cayley graphs and their lexicographically minimal spanning trees. We don't know whether there exists  a Cayley graph $G$ of a finitely generated countable group with intermediate growth and a lexicographically minimal spanning tree $T$ of $G$ such that $\Ibn(T)<\Igr(T)$.

However there are no periodic trees with intermediate growth rate. 
\begin{proposition}\label{prop: periodic tree is either polynomial or exponential}
	Suppose $T$ is an infinite periodic tree. Then either $\br(T)>1$ or there exists an integer $d\geq1$ such that $\big|B(n)\big|=\Theta(n^d)$. Here $\big|B(n)\big|=\Theta(n^d)$ means that the ratio $\big|B(n)\big|/n^d$ is bounded away from zero and infinity. 
\end{proposition}
\begin{proof}
	We give a sketch here and leave the details to interested readers. 
	
	First of all, the periodic tree $T$ is the directed cover of some finite directed graph $G=(V,E)$ based at some vertex $x_0\in V$; see p~82-83 in \cite{LP2016} for a proof of this fact. 
	
	%	With loss of generality we assume for any other vertex $y$ in $V$, there is a directed path from $x_0$ to $y$.
	
	Let $C_1,\ldots,C_m$ be the strongly connected components of $G$ (if for a vertex $v$ there is no directed path from $v$ to itself, then we say $v$ does not belong to any strongly connected component). If there exist some $C_i$ and some $v\in V(C_i)$ such that $v$ has at least two out-going edges in $C_i$, then it is easy to see $\br(T)>1$. 
	Otherwise, each $C_i$ is either a single vertex with a self-loop, or it is a directed cycle. In this case one can prove $\big|B(n)\big|=\Theta(n^d)$ by induction on the size of $V(G)$ (Exercise 3.30 in \cite{LP2016} would be a good warm-up).
	We omit the details of the induction and just point out that in this case
	\[
	d=\max\{ C(\gamma)\colon \gamma \textnormal{ is a self-avoiding directed path in } G \textnormal{ starting from }x_0  \},
	\]
	where $C(\gamma)$ is the number of strongly connected components visited by $\gamma$. 
\end{proof}

\section*{Acknowledgment}

We are grateful to Asaf Nachmias for helpful discussions and Russ Lyons for pointing out an error in a previous draft. We are also grateful to the referees for their careful reading and valuable comments.

 \newcommand\noopsort[1]{}
 \bibliography{subper_ref}
 \bibliographystyle{plain}

\end{document}